\documentclass[10pt,reqno]{amsart}

\usepackage{amsaddr}
\usepackage{etoolbox}
\patchcmd{\section}{\scshape}{\bfseries\scshape}{}{}
\makeatletter
\renewcommand{\@secnumfont}{\bfseries}
\makeatother

\pdfoutput=1 
\usepackage[a4paper,margin=3cm]{geometry}

\usepackage[T1]{fontenc}
\usepackage[utf8]{inputenc}
\usepackage[french,english]{babel}

\usepackage[dvips]{graphicx}
\usepackage{amssymb,amsmath,color}
\usepackage{url}
\usepackage{hyperref}
\usepackage{cancel}
\usepackage{balance}
\usepackage{flushend}

\DeclareMathAlphabet\mathbfcal{OMS}{cmsy}{b}{n}

\newcommand{\BEAS}{\begin{eqnarray*}}
\newcommand{\EEAS}{\end{eqnarray*}}
\newcommand{\BEA}{\begin{eqnarray}}
\newcommand{\EEA}{\end{eqnarray}}
\newcommand{\BEQ}{\begin{equation}}
\newcommand{\EEQ}{\end{equation}}
\newcommand{\BIT}{\begin{itemize}}
\newcommand{\EIT}{\end{itemize}}
\newcommand{\BNUM}{\begin{enumerate}}
\newcommand{\ENUM}{\end{enumerate}}
\newcommand{\BA}{\begin{array}}
\newcommand{\EA}{\end{array}}

\newcommand{\argmin}{\mathop{\rm argmin}}

\newcommand{\Tr}{\mathop{ \rm Tr}}

\newtheorem{theorem}{Theorem}

\newcommand{\e}{\varepsilon}

\newcommand{\Ee}{\mathbb{E}}
\newcommand{\Rr}{\mathbb{R}}
\newcommand{\Ss}{\mathbb{S}}
\newcommand{\Nn}{\mathbb{N}}

\newcommand{\cA}{\mathcal{A}}

\newcommand{\cF}{\mathcal{F}}

\newcommand{\cX}{\mathcal{X}}

\newcommand{\cU}{\mathcal{U}}

\newcommand{\cH}{\mathcal{H}}
\newcommand{\bull}{\raisebox{-0.25ex}{\scalebox{1}{$\cdot$}}}

\allowdisplaybreaks 

\title{Infinite-Dimensional Sums-of-Squares\\ for Optimal Control}

\author{Eloïse Berthier, Justin Carpentier, Alessandro Rudi \and Francis Bach%
}

\address{Inria - Ecole Normale Supérieure\\ 
	PSL Research University, Paris, France}

\email{firstname.lastname@inria.fr}

\begin{document}

\thispagestyle{empty}
\pagestyle{empty}

\begin{abstract}
    We introduce an approximation method to solve an optimal control problem \textit{via} the Lagrange dual of its weak formulation. It is based on a sum-of-squares representation of the Hamiltonian, and extends a previous method from polynomial optimization to the generic case of smooth problems. Such a representation is infinite-dimensional and relies on a particular space of functions -- a reproducing kernel Hilbert space -- chosen to fit the structure of the control problem. After subsampling, it leads to a practical method that amounts to solving a semi-definite program. We illustrate our approach by a numerical application on a simple low-dimensional control problem.
\end{abstract}

\maketitle

\section{Introduction}

The continuous-time optimal control problem (OCP) is a versatile framework modelling a wide variety of nonlinear systems and optimization criteria, with countless industrial applications, including aerospace~\cite{trelat2012optimal} or robotics~\cite{murray2017mathematical}. Developing  efficient numerical methods for solving such a general problem is a daunting task, especially for high-dimensional systems. Among current methods,  \textit{indirect methods}  exploit optimality criteria derived from Pontryagin's maximum principle and give precise results but need to be initialized properly, while the less accurate \textit{direct methods}  reformulate the  problem as a nonlinear program without specific initialization requirements~\cite[Chapter~9]{trelat2005controle}. 

In this paper, we focus on a direct method that computes the optimal value function of the problem as the maximal subsolution of the Hamilton-Jacobi-Bellman (HJB) equation. It is described in~\cite{hernandez1996linear, lasserre2008nonlinear} and is obtained by taking the dual of the \textit{weak formulation} of the OCP,  involving occupation measures~\cite{vinter1993convex}. In~\cite{lasserre2008nonlinear}, the numerical resolution of this formulation is based on polynomial optimization~\cite{lasserre2015introduction}, and hence is restricted to polynomial dynamics and cost functions, with a possibly costly  extension to smooth functions, involving a  hierachy of semi-definite programs (SDPs).

Our main contribution is to extend the numerical method of~\cite{lasserre2008nonlinear} to non-polynomial and smooth OCPs. To this end, we consider a space of smooth functions, called a reproducing kernel Hilbert space (RKHS)~\cite{aronszajn1950theory}, and use representations of non-negative functions in this space~\cite{marteau2020non} with sums-of-squares. This is detailed in Sections~\ref{background} and~\ref{dense}. In Section~\ref{tight}, we notably prove that this representation is exact in two particular cases: the time-invariant linear quadratic regulator (LQR) and smooth control-affine systems. This results in a practical numerical method derived in Section~\ref{formulation},  that requires to solve an SDP to approximate the optimal value function of the OCP. Finally, in Section~\ref{experiment} we illustrate the practical application of this method to a simple two-dimensional OCP.

\section{Background}
\label{background}

First, we introduce the three building blocks that are then combined in Section~\ref{dense} to design our approximation method.

\subsection{Formulation of OCP with Maximal Subsolutions of HJB}

Let $\cX$ and $\cU$ be compact subsets of $\Rr^d$ and $\Rr^p$, for integers $d, p \geq 1$, and assume that $\cU$ is convex. We define the dynamics $f:[0, T] \times \cX \times \cU \rightarrow \Rr^d$, the running cost $L: [0, T] \times\cX \times \cU \rightarrow \Rr$, and the terminal cost $M:\cX \rightarrow \Rr$ occurring at the fixed terminal time $T > 0$. In addition, we assume that the state trajectories remain in the compact set~$\cX$, and that sampling points from $\cX$ is easy. We do not consider problems with explicit state-constraints which are left as future work.

Assume the existence of a smooth function  ${V^* \in  C^{1}(\cX_T)}$, meaning that $V^*$ is smoothly differentiable on \mbox{$\cX_T:=[0, T]\times \cX$}, that is a solution of the HJB equation: $\forall (t, x) \in  \cX_T$,
\begin{align*}
   & \frac{\partial V^*}{\partial t}(t, x) + \inf_{u\in \cU} \left(L(t, x, u) + \nabla V^*(t, x)^\top f(t, x, u) \right) = 0 \\
        &  V^*(T, x) = M(x) , \label{hjb} \tag{HJB}
\end{align*}
where $\nabla V^*$ refers to the gradient of $V^*$ w.r.t.~$x$. For $t \in [0, T]$, let $\cU_t$ the set of admissible controls such that $\forall s \in [t, T], (x(s), u(s)) \in \cX \times \cU$. Then $V^*$ is the value function~\cite{liberzon2011calculus} of the following OCP: $\forall (t_0, x_0) \in \cX_T$,
    \begin{align*}
         & V^*(t_0, x_0) = \inf_{u \in \cU_{t_0}}  \int_{t_0}^T L(t, x(t), u(t)) \text{d}t + M(x(T)) \\ & \forall t\in[t_0, T], ~  \dot x(t) = f(t, x(t), u(t)), \quad x(0) = x_0.
        \tag{OCP} \label{eqn:OCP}
    \end{align*}
Let $\mu_0$ be a probability measure on $\cX$. We are interested in the value of the stochastic initial point problem, where $x_0$ is drawn according to $\mu_0$, that is, $ \Ee_{x_0 \sim \mu_0} \left[ V^*(0, x_0) \right] = \int V^*(0, x_0) \text{d}\mu_0 (x_0).$
    

In this paper, instead of directly looking for solutions of the HJB equation, we will focus on the following alternative problem~(\ref{eqn:P}), namely, finding a maximal subsolution of HJB:
    \begin{align*}
    & \sup_{V \in C^1(\cX_T) } \int V(0, x_0) \text{d}\mu_0(x_0) \\
    \forall (t, x, u),  ~ & \frac{\partial V}{\partial t}(t, x) + L(t, x, u) + \nabla V(t, x)^\top f(t, x, u) \geq 0 \\
        \forall x, ~ & V(T, x) \leq M(x) .
        \tag{P} \label{eqn:P}
    \end{align*}

This is the dual of the \textit{weak formulation} of the OCP with occupation measures, which is a linear program  over the space of measures~\cite{vinter1993convex,lasserre2009moments,kamoutsi2017infinite}. Moreover, subsolutions of HJB also play a key role in the theory of viscosity solutions~\cite{crandall1983viscosity} of partial differential equations. The first constraint in~(\ref{eqn:P})  is the positivity of a certain  Hamiltonian associated to $V$, namely: 
$$H(t, x, u) := \frac{\partial V}{\partial t}(t, x) + L(t, x, u) + \nabla V(t, x)^\top f(t, x, u).$$
If $V = V^*$ is the optimal value function, then the optimal controller $u=u^*(t,x)$ minimizes the positivity constraint and for all $(t, x) \in \cX_T$, $H^*(t, x, u^*(t, x)) =0$.
    
Our goal is to find an approximate solution $V$ of~(\ref{eqn:P}). Under some additional assumptions, (\ref{eqn:P}) is equivalent to the OCP. Such regularity and convexity assumptions were first studied by~\cite{vinter1993convex} and are detailed in~\cite{lasserre2008nonlinear}. In this particular case, the value of problem~(\ref{eqn:P}) coincides with the one of the stochastic initial point problem: $$\sup P = \Ee_{x_0 \sim \mu_0} \left[ V^*(0, x_0) \right].$$

\subsection{Parameterization of the Value Function}
\label{param}

A first difficulty in problem~(\ref{eqn:P}) is searching $V$ in the infinite-dimensional set $C^1(\cX_T)$. One option is to search~$V$ in a finitely-parameterized set $\cF_\Theta$. A  common practice, notably in approximate dynamic programming and  reinforcement learning~\cite[Chapter~9]{sutton2018reinforcement}, is to use a linear approximation of~$V$, with a feature vector $\psi(t, x) \in \Rr^m$ and a parameter~$\theta$ in a convex subset~$\Theta$ of $\Rr^m$, for $m \geq 1$.

Since we assumed that $V^*$ is a solution of the HJB equation, we can restrict that search space in~(\ref{eqn:P}) to functions~$V$ such that $V(T, .) = M(.)$. Then we assume that the parameterization is such that $\forall \theta, V_\theta(T, .) = M(.)$, so that we can remove the explicit constraint in~(\ref{eqn:P}). Hence our parameterized set is:
$$\cF_\Theta:= \{(t, x) \mapsto V_\theta(t, x) = \theta^\top \psi(t, x) + M(x) ~|~ \theta \in \Theta \} ,$$
with $\psi$ such that $\psi(T, .) = 0$. To simplify the evaluations of $\nabla V_\theta$ and $\frac{\partial V_\theta}{\partial t}$, it is convenient (but not necessary) to use a separable feature vector $\psi(t, x) = \kappa(t) \varphi(x)$, with $\kappa(T) = 0$.

\subsection{Representing Non-Negative Functions as Sum-of-Squares}

Problem~(\ref{eqn:P}) is constrained by a dense set of inequalities indexed by $(t, x , u) \in [0, T]\times  \cX \times \cU$, which cannot be directly handled by numerical algorithms. Hence we look for a finite -- possibly approximate -- representation of the non-negative function $(t, x, u) \mapsto H(t, x, u) \geq 0$.

If $f$, $L$, $M$ are polynomials, one way is to use sum-of-squares (SoS) polynomials~\cite{lasserre2015introduction}, \textit{i.e.}, to represent $H(t, x, u)$ as the sum of the squares of polynomials of a given degree. This is a sufficient but not necessary condition for being a non-negative polynomial. This technique has been applied to problem~(\ref{eqn:P}) in~\cite{lasserre2008nonlinear}, although it is presented in its dual version using the method of moments~\cite{henrion2013optimization}. In any case, this representation is not exact in general and gives a lower approximation of~(\ref{eqn:P}). To numerically solve the problem, one needs to build a hierarchy of SDPs obtained by this SoS representation with polynomials of increasing degree~$r$. Under generic conditions, this hierarchy converges to the value of~(\ref{eqn:P}), but the speed of convergence is not explicitly controlled~\cite{robotchapter}. This is a potentially critical issue, since the size of the SDP at rank $r$ is defined by the number of monomials of degree less than $r$ in the dimension of $(t, x, u)$, which is $d+p+1+r \choose r$, a quantity growing exponentially with~$r$.

 In this paper, we opt for another option inspired by recently-introduced machine learning techniques~\cite{rudi2020finding}: representing a non-negative function as a SoS in a reproducing kernel Hilbert space (RKHS). Hereafter, we briefly define an RKHS and mention its main properties, and refer to~\cite{paulsen_raghupathi_2016} for a  thorough description. Consider a set $E$, a function $k : E \times E \rightarrow \Rr$ is a positive definite kernel if $\forall n \in \Nn,  y_1,\dots, y_n \in E$, the matrix ${K := (k(y_i, y_j))_{i,j=1}^n}$ is positive semi-definite. Associated to a positive definite kernel $k$, there exists a unique RKHS $\cH$, a Hilbert space of functions $E \rightarrow \Rr$, with an inner product $\langle \cdot, \cdot \rangle_\cH$, such that the following properties hold (the second one is the so-called ``reproducing property''):
\begin{itemize}
   \item $\forall y \in E$, $k_y := k(y, \cdot) \in \cH$;
   \item $\forall g \in \cH$,  $\forall y\in E$, $\langle g, k_y \rangle_\cH = g(y)$.
\end{itemize}
In addition, there exists a feature map $\Phi:E \rightarrow \cH$, possibly infinite dimensional, defined by $\Phi(y) = k_y$, which maps a point in $E$ to a function in $\cH$, and in particular we have $k(y, y')=\langle \Phi(y), \Phi(y') \rangle_\cH$. Conversely, any feature map defines an RKHS associated to the former kernel.

Here we mention two classical kernels from the RKHS literature, which we use later to build our function representations. Assume that~$E$ is a subset of $\Rr^\ell$, for $\ell \geq 1$. The polynomial kernel of degree~$r$ is defined on $E \times E$ by $k(y, y') = (1 + y^\top y')^r$, and the corresponding embedding~$\Phi(y)$ is the vector of $\ell +r \choose r$ multivariate monomials of degree less than~$r$. In this case~$\cH$ is finite-dimensional. The exponential kernel is defined by $k(y, y') = \exp(- \| y - y' \|_2/\sigma)$, with $\sigma>0$. If $E$ is bounded and has locally Lipschitz boundary, the corresponding RKHS is the Sobolev space of functions whose weak-derivatives up to order $s=\ell/2+1/2$ are square-integrable~\cite{berlinet2011reproducing}. 

Functions that are SoS in an RKHS $\cH$ can be represented using an infinite-dimensional positive semi-definite operator~\cite[Corollary~1]{rudi2020finding}. Indeed, assume that a function $h \in \cH$ is written as a sum-of-squares of functions $h_j \in \cH$: $$\forall y \in E, \quad h(y) = \sum_{j=1}^q h_j(y)^2. $$
For any $v, w \in \cH$, we have: $$\langle w, (v \otimes v) w \rangle = \langle w, v v^* w \rangle = \Tr(w^* v v^* w) = (\langle v, w \rangle)^2,$$
where $\otimes$ denotes the outer product. Because of the reproducing property: $\forall j \in \{1,...,q\}$, $y \in E$,  
 \begin{align*}
     h_j(y)^2 &= (\langle \Phi(y), h_j \rangle_\cH)^2 = \langle \Phi(y), (h_j \otimes h_j) \Phi(y) \rangle . 
 \end{align*}
 And then: $$h(y) = \langle \Phi(y), \cA \Phi(y) \rangle_\cH, $$ with $\cA :=  \sum_{j=1}^q h_j \otimes h_j  \in \Ss_+(\cH)$ and $\cA$ has rank at most~$q$, where $\Ss_+(\cH)$ is the set of bounded self-adjoint positive semi-definite operators on~$\cH$.

In~\cite{marteau2020non}, this SoS representation in an RKHS is used to model non-negative functions, \textit{e.g.}, for signal processing or statistics applications. In some cases, \textit{e.g.}, in Sobolev spaces as we will see hereafter, the representation is exact in the sense that all non-negative functions can be written as a SoS in~$\cH$, whereas the polynomial SoS representation is tight only for a restricted class of polynomials~\cite{lasserre2009moments}. Besides, SoS polynomials are a particular case of what has just been described, if~$k$ is the polynomial kernel. In the rest of the paper, we will extend the method of~\cite{lasserre2008nonlinear} from polynomials to any RKHS, at the expense of possibly infinite-dimensional representations.

\section{Dense Set of Inequality Constraints}
\label{dense}

In this section, we start by providing a basic relaxation of problem~(\ref{eqn:P}) and a motivation for preferring a SoS representation of the non-negativity constraints in~(\ref{eqn:P}). Then we present the resulting problem and its main features.

\subsection{Relaxed formulation by subsampling}
\label{LPmethod}

A straightforward relaxation of~(\ref{eqn:P}) is obtained by finitely subsampling the non-negativity constraints. Let us sample values of $(t^{(i)}, x^{(i)}, u^{(i)})_{i \in I}$ in $[0, T] \times \cX \times \cU$, with $I$ a finite set of cardinality~$n \geq 1$. For simplicity, assume that~$\mu_0$ is the mean of $n_0$ Diracs at points $\{x_0^{(i)}\}_{i \in \{1,\dots,n_0 \}}$. Besides, let us use a linear parameterization of $V$ as described in Section~\ref{param}, with $\Theta=\Rr^m$. We then  obtain a linear program, with a possibly unbounded solution in the overparameterized setting ($m \gg n$). To circumvent this effect, we add a quadratic regularizer on~$\theta$ with parameter $\lambda_\theta \geq 0$, and obtain the following problem:
  \begin{align*}
    & \sup_{\theta \in \Rr^m} \frac{1}{n_0} \sum_{k=1}^{n_0} \theta^\top \psi(0, x^{(k)}_0) + M (x^{(k)}_0) - \lambda_\theta \|\theta\|_2^2 \\
    & \forall i\in I,  \quad  \theta^\top \frac{\partial \psi}{\partial t}(t^{(i)}, x^{(i)}) + L(t^{(i)}, x^{(i)}, u^{(i)}) 
    \tag{LP} \label{eqn:LP}\\& + \left(\theta^\top \nabla \psi(t^{(i)}, x^{(i)}) + \nabla M(x^{(i)})  \right)^\top f(t^{(i)}, x^{(i)}, u^{(i)}) \geq 0 .
    \end{align*}
Although this is not exactly a linear program if $\lambda_\theta>0$, it can still can be solved easily by standard solvers, and we will refer to it as the LP problem. A similar finite-dimensional LP formulation has been proposed in~\cite{gatisgory} for discounted infinite-horizon control problems. It is part of a long series of LP formulations for optimal control (see, \textit{e.g.},~\cite{gaitsgory} and references therein), for dynamic programming and more recently for reinforcement learning~\cite{mehta}.

This problem will be used as a baseline in Section~\ref{experiment}, to be compared with the SoS formulation below. It is a relaxation that gives an upper-bound on~(\ref{eqn:P}), but it is not straightforward to relate the number of samples~$n$ and how they are spread to the quality of the approximation. Yet in the example below, this can be evaluated explicitly.
    
\textit{Example~1:} Let $g : \Rr^p \rightarrow \Rr$ be a smooth function  with a unique minimizer $u^*$. With \mbox{$T=1$, $L:(t, x, u)\mapsto g(u)$}, ${f=0}$ and $M=0$,  then $V^*(t, x) = (1-t) g(u^*)$ and solving the OCP is essentially equivalent to finding the global minimizer of~$g$.  If~$V$ is parameterized by $V_\theta(t, x) = \theta (1-t)$, the LP formulation with $\lambda_\theta=0$ writes:
    $$\sup \theta \text{ s.t. } \forall i \in I, -\theta + g(u_i) \geq 0, $$
    which is readily solved by $\theta=\min \{g(u_i)\}_{i \in I}$. In general, this method requires $O(\e^{-p})$ samples to approximate $g(u^*)$ with precision~$\e$~\cite{novak2006deterministic}, and yet when~$g$ is smooth, this is not an optimal way to perform zero-th order optimization. Indeed,~\cite{rudi2020finding} use a SoS representation of $g-\theta$ to solve this exact problem, and alleviate the curse of dimensionality when~$g$ is sufficiently smooth: the number of samples reduces to $O(\e^{-p/s})$ for $g\in C^s(\Rr^p)$.  In the rest of this paper, we propose to use the same approach and to generalize it to any OCP. We expect similar benefits when~$H$ is smooth.

\subsection{Strengthened formulation by SoS representation}
\label{SOS}

Consider an RKHS $\cH$ of real-valued functions on $[0, T]\times \cX \times \cU$, with positive definite kernel $k$, and $\Phi:E\rightarrow \cH$ the corresponding embedding. We use a SoS representation in~$\cH$, or ``kernel SoS'', for the constraint $H \geq 0$ in~(\ref{eqn:P}):
  \begin{align*}
    \sup_{\substack{V \in C^1(\cX_T),\\ \cA \in \Ss_+(\cH) }}~ \int V(0, x_0) \text{d}\mu_0(x_0)  \qquad\qquad \\
    \forall (t, x, u),  ~  H(t, x, u) = \langle \Phi(t, x, u), \cA \Phi(t, x, u) \rangle. \tag{KSOS} \label{eqn:KSOS}
    \end{align*}
    
    
This is a strengthening   of the constraint in~(\ref{eqn:P}), since being SoS is stronger than being non-negative. So in general,~(\ref{eqn:KSOS}) is a lower-approximation of~(\ref{eqn:P}). However, in certain cases,~(\ref{eqn:KSOS}) can be equivalent to~(\ref{eqn:P}), as we will prove in Section~\ref{tight}. A sufficient condition is  the existence of $\cA \in \Ss_+(\cH)$ such that, at the optimal~$V^*$:
$$\forall (t, x, u), ~H^*(t, x, u) = \langle \Phi(t, x, u), \cA \Phi(t, x, u) \rangle. $$

\section{Tight Sum-of-Squares Representations}
\label{tight}

We study the tightness of problem~(\ref{eqn:KSOS}) in two particular cases: the time-invariant LQR and smooth value functions.

\subsection{Case 1: Infinite-Horizon Time-Invariant LQR}
\label{lqr}
    
First we look at a very simple OCP, where every quantity can be computed almost in closed form, and with infinite-horizon so that there is no dependence in~$t$. Let $f(x, u) = A_0x + B_0u$, for $A_0 \in \Rr^{d\times d}$, $B_0\in \Rr^{d \times p}$, with $(A_0, B_0)$ controllable,  $L(x, u) = x^\top Q_0 x + u^\top R_0 u$,  $Q_0 \in \Ss_+(\Rr^{d\times d})$, $R_0 \in \Ss_{+}(\Rr^{p\times p})$, $R_0 \succ 0$. The optimal value function is ${V^*(x) = x^\top S_0 x}$, where $S_0$ is the unique positive semi-definite solution of the algebraic Riccati equation:
$$0 = - Q_0 -A_0^\top S_0 - S_0 A_0 + S_0 B_0R_0^{-1}B_0^\top S_0. $$ The optimal controller is $u^*(x) = - R_0^{-1}B_0^\top S_0 x =: - K_0x$.
\begin{align*}
    &H^*(x, u) = x^\top Q_0 x + u^\top R_0 u + 2 x^\top S_0 (A_0x + B_0u) \\
    &\qquad = u^\top R_0 u + x^\top S_0B_0u + u^\top B_0^\top S_0 x  + x^\top S_0B_0K_0 x.
\end{align*}
This is a SoS of degree-one polynomials in $(x, u)$: 
\begin{align*}
    H^*(x, u) &= (u + K_0 x)^\top R_0 (u+K_0 x) \\
    &= \begin{pmatrix} x^\top & u^\top \end{pmatrix} \begin{pmatrix} K_0^\top \\ I_p \end{pmatrix} R_0 \begin{pmatrix} K_0  &  I_p \end{pmatrix} \begin{pmatrix} x\\u \end{pmatrix} \\[-.2cm]
    &=  \sum_{j=1}^p [q_j(x, u)]^2,
\end{align*}
with $q_j(x, u) := [R_0^{1/2}]_{j\bull} \begin{pmatrix} K_0 &  I_p \end{pmatrix} \begin{pmatrix} x \\ u \end{pmatrix}$.

Hence, an infinite-horizon, time-invariant LQR with unknown parameters can be equivalently expressed by:
\begin{align*}
     \sup_{V \in C^1(\cX),~ N \succeq 0 } ~\int V(x_0) \text{d}\mu_0(x_0) \qquad\qquad\qquad \\
    \forall (x, u),  \quad  
    L(x, u) + \nabla V(x)^\top f(x, u) = \begin{pmatrix} x \\ u \end{pmatrix}^\top N \begin{pmatrix} x \\ u \end{pmatrix}.
    \end{align*}

In the next section, we prove that similar SoS constructions exist for sufficiently smooth OCPs, possibly with an infinite-dimensional embedding (\textit{v.s.} a $(d+p)$-dimensional one here).

\subsection{Sum-of-Squares Decomposition with Smooth Functions}

We show that, for smooth and control-affine OCPs, $H^*$ is a SoS of smooth functions. Let $\Omega_2 := \textnormal{Int}(\cU)$, $\Omega_1 := \textnormal{Int} \{(t, x) \in \cX_T ~|~ \argmin_{u \in \cU} H^*(t, x, u) \subset \Omega_2 \}$, and $\Omega := \Omega_1 \times \Omega_2$.

\begin{theorem} 
Let $s \in \Nn$, $s \geq 1$. Assume that:

\noindent $\bullet$ $f$ is control-affine: $\forall (t, x, u) \in [0, T] \times \cX \times \cU$, $$f(t, x, u) = g(t, x) + B(t, x)u.$$
$\bullet$ For all $(t, x) \in \Omega_1$, $u \mapsto L(t, x, u)$ is twice differentiable on $\Omega_2$ and strongly convex:  $\nabla^2_u L(t, x, u) \succcurlyeq \rho I$ for some $\rho >0$, and $(t, x, u) \mapsto \nabla^2_u L(t, x, u) \in C^s(\Omega)$.

\noindent $\bullet$ $(t, x, u) \mapsto \nabla_u L(t, x, u) + B(t, x)^\top \nabla_x V^*(t, x) \in C^s(\Omega)$.

\noindent Then there exist~$p$ functions $(w_j)_{1\leq j \leq p} \in C^s(\Omega)$ such that:
$$\forall (t, x, u) \in \Omega, \quad H^*(t, x, u) = \sum_{j=1}^p w_j(t, x, u)^2. $$
\end{theorem}

The proof is in the appendix. This result motivates the use of exponential kernels, inducing a Sobolev space RKHS, to represent the non-negativity constraints in~(\ref{eqn:P}) for smooth OCPs. When $s> d/2+1$, by applying a technique similar to the one used under Corollary 2 of~\cite{rudi2020finding}, it is possible to obtain a SoS representation in terms of the exponential kernel. Then~(\ref{eqn:KSOS}) is equivalent to~(\ref{eqn:P}).

\subsection{Stochastic Smoothing of the Optimal Value Function}

However in general, $V^*$ is not necessarily smooth (e.g. minimal time problems), nor is~$u^*$ (e.g. \textit{bang-bang} controllers that are not even continuous). Here, we provide a generic technique to give some regularity to~$V^*$. For $\eta > 0$, we can define a perturbed version of the control system~\cite{fleming2012deterministic}, where the state is a random variable $X_t$:
$$ \text{d}X_t = f(t, X_t, u_t)\text{d}t + (2\eta)^{1/2} \text{d}B_t, $$
where $B_t$ is a standard Brownian motion independent of $X_t$. We define the optimal value function, with $X_0 = x_0$:
$$V^\eta(t_0, x_0) = \inf_{u\in \cU} \Ee \left[ \int_{t_0}^T L(t, X_t, u_t) \text{d}t + M(X_T) \right]. $$
$V^\eta$ is the unique $C^{1,2}(\cX_T)$  ($C^1$ in $t$, $C^2$ in $x$)  solution~\cite{fleming2012deterministic} of the following regularized HJB equation: $\forall (t, x) \in \cX_T$,
\begin{align*}
      \inf_{u\in \cU} ~ \{ L(t,x,u)  ~+ & \nabla V^\eta(t, x)^\top f(t, x, u) \} \\ ~+   & \frac{\partial V^\eta}{\partial t}(t, x) + \eta \Delta V^\eta(t, x)= 0,
\end{align*}
with $V^\eta(T, x) = M(x)$. $\Delta V^\eta$ refers to the Laplacian with respect to $x$. Contrary to HJB, the solutions are at least $C^{1,2}(\cX_T)$ because this is a quasilinear parabolic partial differential equation~\cite{lieberman1996second}. The regularization $\eta \Delta V$ is a vanishing viscosity term, and the optimal controller is still in $\argmin_u L+ \nabla {V^\eta}^\top f$.  Generically, $V^\eta$ converges to $V^*$ as $\eta \rightarrow 0$, following a reasoning similar to the theory of viscosity solutions~\cite{crandall1983viscosity}.

\section{SDP Formulation and its Numerical Resolution}
\label{formulation}

\subsection{Finite-Dimensional Formulation via Subsampling}

Similarly to the~(\ref{eqn:LP}) formulation that relaxes~(\ref{eqn:P}), we will now derive a relaxation of problem~(\ref{eqn:KSOS}), which is another relaxation of problem~(\ref{eqn:P}) if the SoS representation of~$H^*$ is tight. Going through~(\ref{eqn:KSOS}) as an intermediate step will help to exploit the structure of~(\ref{eqn:P}). Using a parameterization of~$V$ in~$\cF_\Theta$ with $\Theta=\Rr^m$, and a set of sampled points $(t^{(i)}, x^{(i)}, u^{(i)})_{i \in I}$  in $[0, T] \times \cX \times \cU$, with $|I|=n$, we obtain:
   \BEAS
  \!\!\!&\sup_{\cA \in \Ss_+(\cH), \theta\in \Theta}  & c^\top \theta  - \lambda_\theta \| \theta\|_2^2 - \lambda \Tr (\cA) + C\\ 
 \!\!\! &\!\!\!\mbox{ such that }   &  \forall i\in \{1,\dots,n\} , \\ && \hspace{-2cm} b_i + a_i^\top \theta =  \langle \Phi(t^{(i)}, x^{(i)}, u^{(i)}) , \cA \Phi(t^{(i)}, x^{(i)}, u^{(i)}) \rangle,
\EEAS
with $c:=\sum_{i=1}^n \mu_0^{(i)}   \psi(0, x^{(i)})$, $C := \sum_i \mu_0^{(i)} M(x^{(i)})$, $b_i := L(t^{(i)}, x^{(i)}, u^{(i)}) + \nabla M(x^{(i)})^\top  f(t^{(i)}, x^{(i)}, u^{(i)}) + \eta \Delta M(x^{(i)})$, and $a_i := J_\psi(t^{(i)}, x^{(i)}) f(t^{(i)}, x^{(i)}, u^{(i)}) + \frac{\partial \psi}{\partial t}(t^{(i)}, x^{(i)}) + \eta \Delta \psi (t^{(i)}, x^{(i)})$, where $J_\psi$ denotes the Jacobian matrix of $\psi$ with respect to~$x$ only. Note that we integrate the stochastic smoothing process in this formulation, with parameter~$\eta$ that can be eventually set to 0.

The regularization parameter $\lambda>0$ controls the trace of the infinite-dimensional operator~$\cA$, and allows for subsampling to provably recover the non-subsampled program when~$n$ tends to infinity, and~$\lambda$ goes to zero at the proper rate (see \cite{rudi2020finding} for the precise dependence). In the limit $\lambda \rightarrow 0$, we recover the LP formulation where we assume nothing about the SoS representation of~$H^*$ in~$\cH$. 

Both the operator $\cA$ and the $\Phi(t^{(i)}, x^{(i)}, u^{(i)})$ can be infinite dimensional, depending on the RKHS $\cH$. Yet, following~\cite{rudi2020finding}, we can  reformulate the problem equivalently in finite dimension. Using the representer theorem in~\cite{marteau2020non}, one can prove that~$\cA$ can be sought in the form: for $D\in \Rr^{n\times n}$, $D\succeq 0$, $$ \cA = \sum_{i,j=1}^n D_{ij} \Phi(t^{(i)}, x^{(i)}, u^{(i)}) \otimes \Phi(t^{(j)}, x^{(j)}, u^{(j)}).$$
Simple computations detailed in~\cite{rudi2020finding} show that:
$$
 \left\{
    \begin{array}{l}
            \forall i, ~  \langle \Phi(t^{(i)}, x^{(i)}, u^{(i)}) , \cA \Phi(t^{(i)}, x^{(i)}, u^{(i)}) \rangle = \left[KDK\right]_{ii}, \\
       \Tr(\cA) = \Tr(DK),
    \end{array}
\right.
$$ where $K$ is the kernel matrix  with entry $(i, j)$ equal to $k\left( (t^{(i)}, x^{(i)}, u^{(i)}), (t^{(j)}, x^{(j)}, u^{(j)}) \right)$. Assume that $K\succ 0$. We denote by $K=R^\top R$ the Cholesky decomposition of~$K$, with~$R$ an invertible upper-triangular matrix.

Let $B:=RDR^\top$ and for $1\leq i \leq n$, $\Phi_i := R_{\bull i}$. Then:
$$
 \left\{
    \begin{array}{l}
       \Tr(B) = \Tr(DK) = \Tr(\cA), \\
       \left[KDK\right]_{ii} = \left[R^\top B R\right]_{ii} = \Phi_i^\top B \Phi_i.
    \end{array}
\right.
$$ The problem can now be reformulated as a finite-dimensional SDP over the positive semi-definite matrix~$B \in \Rr^{n\times n}$:
\BEAS
  \!\!\!&\sup_{B \succcurlyeq 0, \theta\in \Rr^m}  & c^\top \theta - \lambda_\theta \| \theta\|_2^2 - \lambda \Tr (B) + C   \label{sdp}  ~~~~~~\text{(SDP)} \label{eqn:SDP}\\ 
 \!\!\! &\!\!\!\mbox{ such that }   &  \forall i\in \{1,\dots,n\} , \ b_i + a_i^\top \theta  = (\Phi_i )^\top B \Phi_i .
\EEAS
An important question is to estimate the  number of subsampled inequalities sufficient  to ensure that (SDP)$\ \Leftrightarrow$~(\ref{eqn:KSOS}). If nothing is assumed on the structure of $H$, as in the LP method, this number is infinite. In contrast, the kernel SoS representation can reduce it or make it finite. If $H$ is a polynomial of degree $2r$, $k$ is the polynomial kernel of degree $r$, then $n \geq 2r$ distinct sampled points are enough to interpolate~$H$, and~(SDP)$\	\Leftrightarrow$~(\ref{eqn:KSOS}). Another example is global optimization of smooth functions (see \textit{Example~1}) with the exponential kernel. We refer to~\cite{rudi2020finding} for the analysis of the convergence rates, with a lower dependence in the dimension for the kernel SoS when compared to direct inequality subsampling (corresponding to the LP approach).

\subsection{Interior Point Method with the Damped Newton Method}

Problem~(SDP) can be readily solved by any off-the-shelf SDP solver. However, for large~$n$, this quickly becomes too computationally demanding. Here, we propose a numerical scheme based on~\cite{rudi2020finding} that scales better with the number of subsamples~$n$. First, we introduce a slack variable~$\delta \in \Rr^n$ allowing the constraints to be slightly violated (\textit{e.g.} because $\cF_\Theta$ is not a perfect model), controlled by a large parameter~$\gamma >0$. Second, we introduce a \mbox{log-barrier} term controlled by a small~$\e>0$, useful to form the dual of the SDP. We obtain the following problem:
   \begin{align*}
        \sup_{\substack{B \succcurlyeq 0, \\  \theta, \delta }} ~~ & c^\top \theta - \lambda \Tr (B) - \lambda_\theta \| \theta\|_2^2 - \gamma \| \delta \|^2+ \e \log \det B+C\\ 
\! \mbox{ such}& \mbox{ that }  ~~   \forall i\in \{1,\dots,n\} , \ b_i + a_i^\top \theta  = (\Phi_i )^\top B \Phi_i + \delta_i.
   \end{align*}
The Lagrange dual of this problem reads:
\BEAS
  \!\!\!&\inf_{\alpha\in \Rr^n}  & \sum_{i=1}^n \alpha_i b_i  + \frac{1}{4\lambda_\theta} \sum_{j=1}^m \left( c_j + \sum_{i=1}^n \alpha_i a_{ij} \right)^2 \\ & & - \e \log \det U(\alpha) + \frac{1}{4 \gamma} \| \alpha \|_2^2 + \e n \log (\e/e) + C,
\EEAS
where $U(\alpha) := \lambda I_n + \Phi^\top \text{Diag}(\alpha) \Phi$, and $\Phi:=R^\top$ is the matrix with rows $(\Phi_i)_{1\leq i\leq n}$. Let us call the objective~$F(\alpha)$.

Since $F/\e$ is self-concordant~\cite{boyd2004convex} like in~\cite{rudi2020finding}, we propose to use damped Newton iterations~\cite{nemirovski2004interior} on $F/\e$:
$$\alpha \leftarrow \alpha + \frac{1}{1+\lambda(\alpha)} \Delta \alpha, $$
where $\Delta(\alpha) := - [F''(\alpha)]^{-1} F'(\alpha)$ is the Newton direction and $\lambda(\alpha) := \sqrt{\Delta \alpha^\top F''(\alpha) \Delta \alpha / \e}$ is the Newton decrement. The gradient and Hessian of $F$ are computed by:
\begin{align*}
    \frac{\partial F}{\partial \alpha_i} &= b_i + \frac{1}{2 \lambda_\theta} \sum_{j=1}^m a_{ij}  \left(c_j + \sum_{k=1}^n a_{kj} \alpha_k  \right) + \frac{1}{2\gamma} \alpha_i \\ & \qquad\qquad\qquad\qquad\qquad\qquad  - \e~ \Phi_i^\top U(\alpha)^{-1} \Phi_i .\\
    \frac{\partial^2 F}{\partial \alpha_i \partial \alpha_j} &= \frac{1}{2\lambda_\theta} \sum_{k=1}^n a_{ik} a_{jk} + \e \left[ \Phi_i^\top U(\alpha)^{-1} \Phi_j \right]^2 + \frac{\mathbf{1}_{i=j}}{2 \gamma}.
\end{align*}

At optimum, the value function is recovered by $$\theta^\star = \frac{1}{2 \lambda_\theta}\left(\sum_{i=1}^n \alpha_i^\star a_i  + c \right),$$ and the dual variable~$\alpha^\star$ plays a role similar to an occupation measure~\cite{vinter1993convex}, although it is not necessarily non-negative. To improve numerical stability in the experiments hereafter, we used an homotopy heuristics that progressively decreases the parameters $\lambda_\theta$ and $\e$. Moreover, parallel implementations are possible because no singular value decomposition is needed, only matrix operations and system inversions.

\section{Numerical Example}
\label{experiment}

In this section, we apply the kernel SoS method along with the basic LP method, on a two-dimensional control problem, namely the double integrator with finite horizon.

\paragraph{Setting} The problem is an LQR, as in Section~\ref{lqr}, but with finite-horizon~$T=1$, $d=2$, $p=1$,  $M(x) = \|x\|_2^2 $,
$$A_0 = \begin{pmatrix} 0 & 1 \\ 0 & 0\end{pmatrix}, ~ B_0 = \begin{pmatrix} 0 \\ 1\end{pmatrix}, ~ Q_0 = I_2, ~ R_0 = 0.1~. $$
The optimal value function and controller are $V^*(t, x) = x^\top S(t) x$, $u^*(t, x)=- R_0^{-1}B_0^\top S(t)x =:-K(t) x$, where $S(.)$ is the positive semi-definite solution of $S(T) = I_2$ and:
$$\dot S(t) = - Q_0 -A_0^\top S(t) - S(t) A_0 + S(t) B_0R_0^{-1}B_0^\top S(t). $$

\paragraph{Parameterization of $V$} Let  $V_\theta(t, x) = \theta^\top \psi(t, x)$, where each entry of $\psi$ is a product of basis functions on~$\cX$ and $[0, T]$. Let $\varphi(x) := (1, x_1, x_2, x_1x_2, x_1^2, x_2^2)^\top$, because we know~$V^*$ is quadratic in $x$. For $\kappa$ on $[0,T]$, we only know that it is a smooth function, so we use an approximate basis of the Sobolev space of functions with squared integrable derivatives: a sequence of sines and cosines  with decreasing periods beginning with $2T$ to avoid constraining $V(0, .) = V(T, .)$, and  $\kappa(T) = 0$, ensures that $V(T, .) = M(.) $:
$$\kappa(t) := \left( \frac{1}{\omega} \sin \left(\frac{\omega \pi}{2}\frac{t-T}{T} \right) \right)_{1\leq \omega \leq m_t}^\top ~.$$ 
 Finally, $\psi_{i+6j}(t, x) := \varphi_i(x) \kappa_j(t)$, and $\theta\in \Rr^m$, $m = 6 m_t$. We choose $m_t=10$, for which the performance of the policy of the projection of $V^*$ on $\cF_\Theta$ is almost perfect.

\paragraph{Evaluation} We give two criteria to evaluate the quality of an approximation~$V$. First, the distance to~$V^*$: $\| \bar V - \bar V^*\|^2$, where $\bar V$ is the vector of its values on a regular grid on $[0, T] \times [-1, 1]^2$ with $10\times10\times10$ points. Second, the cost of the policy on a $10\times10$ regular grid of initial points.

\paragraph{Sampling} The set of samples $(t^{(i)}, x^{(i)}, u^{(i)})_{i \in I}$ is built as follows. The $x^{(i)}$ are $n_x$ points in $[-1, 1]^2$ generated by the Sobol sequence~\cite{sobol1967distribution}, the $(u^{(i)})_{1\leq i \leq n_u}$ are on a uniform grid on $[-10, 10]$ and  the $(t^{(i)})_{1\leq i \leq n_t}$ on $[0, T]$. The sample set is the Cartesian product of the three previous ones, and has $n=n_t n_x n_u$ elements. We also use the same samples as initial points $(t_0^{(i)}, x_0^{(i)})$ in the objective function of problem~(\ref{eqn:P}). Note that we have replaced it with $\sum_{i=1}^n V(t^{(i)}, x^{(i)})/n ,$ as we found  it more efficient in our experiments to optimize over~$V$ at intermediate time steps rather than at~$t_0$ only. Indeed, we ultimately evaluate our approximation by the accuracy of~$V$ on the whole $\cX_T$ and not only on $\{t_0\} \times \cX$. In a discrete states and actions setting, this effect is analyzed in~\cite{de2003linear}, where~$\mu_0$ is denoted as ``state-relevance weights''.

\paragraph{Methods} We compare three methods: the LP, the guided SoS and the kernel SoS. The LP method is detailed in Section~\ref{LPmethod}, and as for the kernel SoS method, we add a slackness parameter on the constraints, with a penalization controlled by~$\gamma>0$ ($\gamma \rightarrow \infty$ recovers the original LP). 

The guided SoS method is the same as the~(SDP) problem, except that the embeddings $\Phi_i$ of the samples are replaced by vectors~$\Psi_i$ of fixed dimension, which are computed explicitly, without a kernel. Motivated by the fact that:
$$H^*(t, x, u) = (u + K(t)x)^\top R (u + K(t)x), $$
we choose the embedding vectors as follows:
$$\Psi_i := ( u^{(i)}/10, x^{(i)}, \left( 1/\omega \sin(\omega\pi/2 (t/T-1)) x \right)_{1\leq \omega \leq q_t} )^\top ,$$
where the last~$q_t$ scalar terms (without the vector~$x$) approximately model $K(.)$ as a smooth function of~$t$. For computational efficiency, we choose~$q_t=5$ and we checked that this basis can approximate the entries of~$K(t)$ well. Then we solve an SDP of size $(p+q_t d) \times (p + q_t d)=11\times 11$ instead of $n \times n$ for the kernel version.

The kernel SoS method is as described in the previous sections, with the following kernel:
$$k((t,x,u), (t',x',u')) =  \langle u, u'\rangle/100 + \langle x, x'\rangle \times \exp(- |t-t'|).
 $$
This kernel is also designed to match the shape of $H^*$, with a smooth term in $t$ modelled by the exponential kernel. The matrix~$K$ can be singular, so we replace it by $K+10^{-8} I_n$.

\begin{figure}[t]
    \centering
    \includegraphics[width=10.5cm]{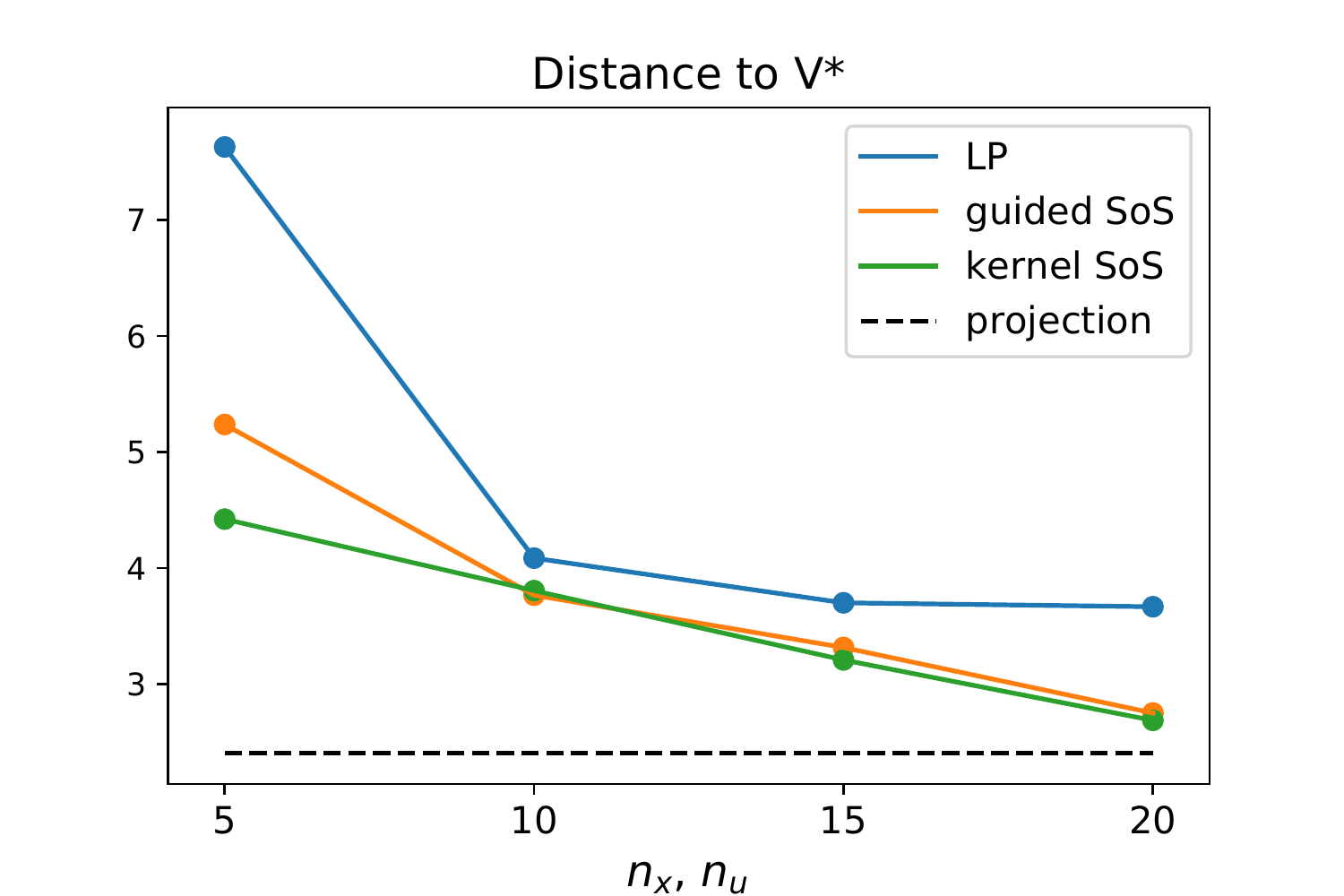}
    \includegraphics[width=10.5cm]{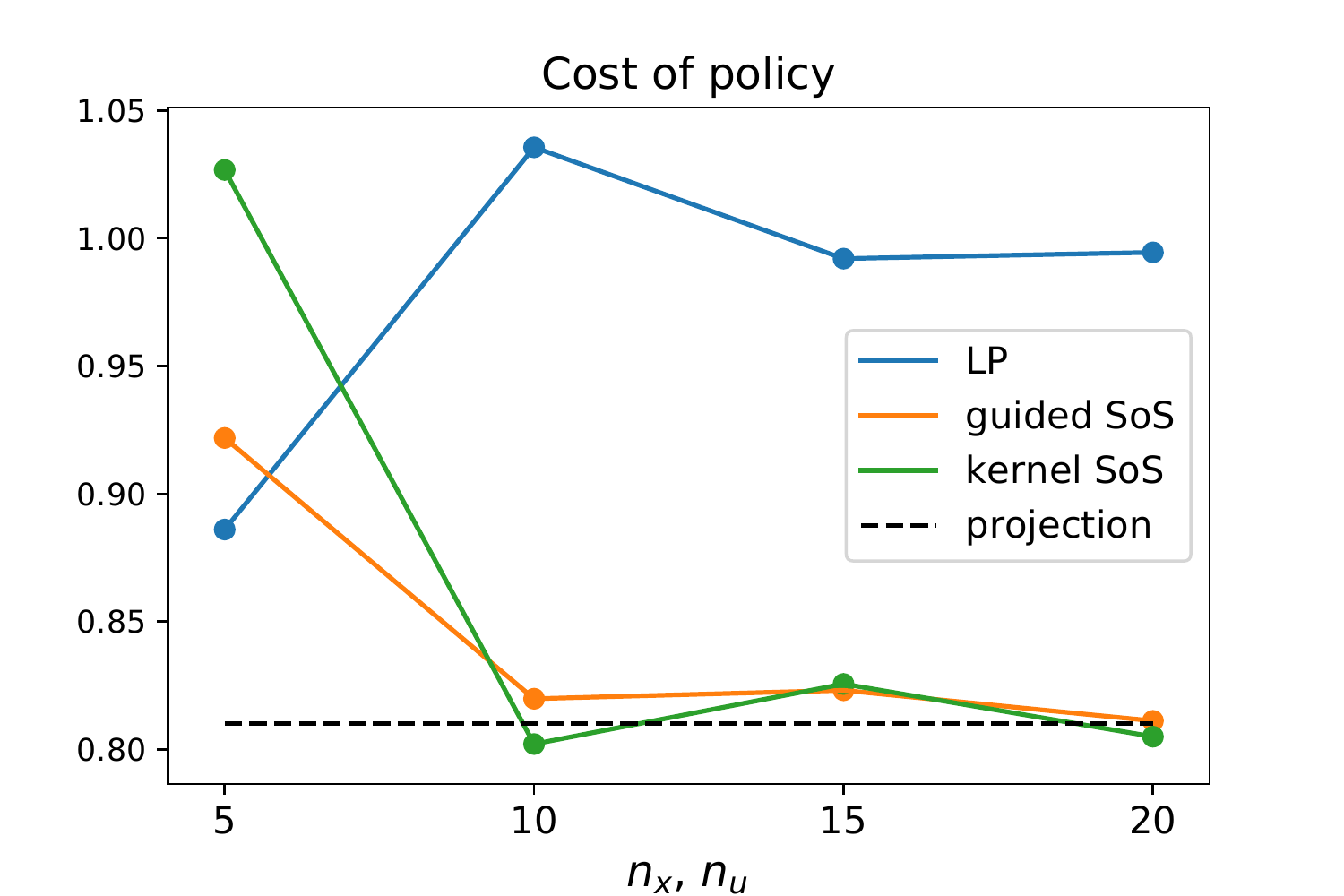}
    \caption{Comparison of the performances of the value function and the policy of the three methods, as a function of the number of samples $n_x=n_u$.}
    \label{fig:results}
\end{figure}

\paragraph{Results} We compare the performance of the three methods to a baseline: the projection of $V^*$ on~$\cF_\Theta$, which is a proxy for the best performance to expect with a fixed~$\cF_\Theta$. We keep the best set of hyper-parameters after a grid search on ($\lambda_\theta$, $\gamma$) for the LP method, and on ($\lambda_\theta$, $\gamma$, $\lambda$) for the two others, with $\e=10^{-4}$. We keep $\eta=0$ and $n_t=20$ in all the experiments, and a varying number $n_x=n_u$ of sample points. For example, with $n_x=n_u=20$, the dual variable~$\alpha$ of the largest problem here has dimension $n=8000$, and solving the numerical problem written in Python takes a few minutes on a standard laptop.

The results are presented in Figure~\ref{fig:results}. The guided and kernel SoS methods perform similar, and better than the LP: they better exploit a fixed number of samples than the LP. Note that the kernel SoS tends to the LP when $\lambda$ tends to 0, hence using a  positive $\lambda$ improves the results. We believe that the design of a kernel adapted to prior knowledge on the problem is crucial to benefit from this effect. Finally, the kernel SoS has the same performance as the guided SoS, but it is computationally more expensive as soon as $n > 11$. Yet the kernel version extends way beyond such fixed finite-dimensional embeddings, to infinite-dimensional embeddings represented by any positive definite kernel, including the exponential kernel, the polynomial kernel and many others.

\section*{Conclusion}

The kernel SoS approximation method  generalizes the polynomial SoS method for OCPs. Like the simple LP method, it is \textit{black-box} in the sense that it is based only on function  evaluations of the dynamics and loss, without requiring any gradients. Moreover, it enables to exploit prior knowledge on the structure of an OCP, by choosing an appropriate kernel. The problem reduces to an SDP, whose size can be computationally limiting, but parallel implementations are possible. There are several sources of approximation in this method: the parameterization $\cF_\Theta$ of $V$ might not be exact, the SoS representation of $H^*$ is not exact in general (although we have proved it is in a few particular cases), and we subsample a finite number of constraints (further work is needed to evaluate the effect of this step). In particular, it seems essential to assess in which cases the   subsampling step is tight with a finite number of samples, or approximately so, in such a way that the overall process gives a certified lower-bound on the OCP, similarly to~\cite{lasserre2008nonlinear}. For all these reasons, the method will probably not reach high precision solutions, but can be used to initialize direct shooting methods, and returns an approximate solution even with very few samples. Furthermore, we believe it is possible to extend the method to also account for state constraints, similarly to~\cite{lasserre2008nonlinear}. One could also parameterize the value function directly in an RKHS. Another interesting extension is to apply the method to Markov decision processes, where we could deal with states or actions that are more complex objects (graphs, trajectories, DNA sequences...) with appropriate kernels.

\section*{Acknowledgements}
  
  This work was supported by the Direction G\'{e}n\'{e}rale de l'Armement, and by the French government under management of Agence Nationale de la Recherche as part of the ``Investissements d'avenir'' program, reference ANR-19-P3IA-0001 (PRAIRIE 3IA Institute). We also acknowledge support from the European Research Council (grants SEQUOIA 724063 and REAL 947908). We thank Didier Henrion for interesting discussions related to this work.

\bibliographystyle{amsalpha}

\bibliography{control}

\appendix 
\section{Proof of Theorem 1}

\begin{proof}
 Consider the Hamiltonian with $f$  control-affine:
\begin{align*}
H^*(t, x, u) = &  ~ \nabla V^*(t, x)^\top g(t, x) + \frac{\partial V^*}{\partial t}(t, x)  \\ & \qquad + L (t, x, u) + \nabla V^*(t, x)^\top B(t, x)u.
\end{align*}
Pontryagin's maximum principle states that:
$$ \forall (t, x) \in \Omega_1, \quad \inf_{u \in \cU}~ H^*(t, x, u) = 0 .$$
Since $L$ is strongly convex in $u$, the minimizer is unique and we call it $u^*(t, x)$. By definitions of $\Omega_1$ and $\Omega_2$, we have a mapping $u^* : \Omega_1 \rightarrow \Omega_2$ and it is characterized by:
\begin{align*}
    \nabla_u H^*(t, x, u) &= \nabla_u L (t, x, u^*(t, x)) + B(t, x)^\top \nabla_x V^*(t, x) \\
    &= 0. 
\end{align*}
Since $(t, x, u) \mapsto \nabla^2_u L(t, x, u)$ is continuous on $\Omega$ and invertible, and $(t, x, u) \mapsto \nabla_u H^*(t, x, u) \in C^s(\Omega)$, then the implicit function theorem ensures that $u^* \in C^s(\Omega_1)$ (see~\cite{schwartz}, Chapter 8, Theorems 25 \& 31).

For $(t, x, u) \in \Omega$, we use Taylor's formula around $u^*(t, x)$:
\begin{align*}
H^*(t, x, u) &= H^*(t, x, u^*(t, x)) \\ & \quad + \nabla_u H^*(t, x, u^*(t, x))^\top (u - u^*(t, x)) \\& \qquad + (u - u^*(t, x))^\top R(t, x, u) (u - u^*(t, x)),
\end{align*}
$$R(t, x, u) := \int_0^1 (1-\tau) \nabla^2_u H^*(t, x, (1-\tau)u^*(t, x) + \tau u) \text{d}\tau.$$
Since $H^*(t,x, u^*(t, x)) = 0$, $\nabla_u H^*(t,x, u^*(t, x)) = 0$ (by definition), and $\nabla_u^2 H^*(t, x, \cdot) = \nabla^2_u L(t, x, \cdot)$, we have:
\begin{align*}
H^*(t, x, u) =  (u - u^*(t, x))^\top R(t, x, u) (u - u^*(t, x)), \text{~and}
\end{align*}
$$R(t, x, u) = \int_0^1 (1-\tau) \nabla^2_u L(t, x, (1-\tau)u^*(t, x) + \tau u) \text{d}\tau \succcurlyeq \frac{\rho}{2} I. $$
For $(t, x, u) \in \Omega$, $R(t, x, u)$ has a positive-definite square root $\sqrt{R(t, x, u)}$. Also, $\forall \tau \in[0, 1]$, $(1-\tau)u^*(t, x) + \tau u \in \Omega_2$ because $\textnormal{Int}(\cU)$ is convex like $\cU$.

Since $\forall i, j, \frac{\partial^2 L}{\partial u_i \partial u_j} \in C^s(\Omega)$, $u^* \in C^s(\Omega_1)$, and $\sqrt{\cdot}$ is $C^\infty$ on $\{M~|~M^\top = M, M\succcurlyeq \frac{\rho}{2} I \}$, then $r_{i,j} : (t, x, u) \mapsto e_i^\top \sqrt{R(t, x, u)} e_j \in C^s(\Omega)$, and we have the decomposition:
$$H^*(t, x, u) = \sum_{i=1}^p w_i(t, x, u)^2, \text{~~ with} $$
\begin{align*}
    w_i(t, x, u) &:= \sqrt{R(t, x, u)}_{i\bull} (u - u^*(t, x))  \\ &= \sum_{j=1}^p r_{i,j}(t, x, u) \left(e_j^\top (u - u^*(t, x)) \right),
\end{align*}
and each $w_i \in C^s(\Omega)$.
\end{proof}


  \end{document}